# CORRECTION BOUNDS ON MEASURES SATISFYING MOMENT CONDITIONS

By Jean B. Lasserre

In the statement of Theorem 4.1, page 1130 in Lasserre (2002), one should replace "*as* $r \to \infty$," by "*as* $r \to \infty$, *and provided* $\varepsilon(r) \downarrow 0$ *sufficiently slowly*."

Indeed, an arbitrary nonincreasing sequence $\varepsilon(r) \downarrow 0$ is defined on page 1130. Then, in the proof of Theorem 4.1, pages 1134–1136 in Lasserre (2002), with $\varepsilon > 0$ fixed, arbitrary, one obtains the identity (5.7), after which one defines $r := \max[r_0, r_1, r_2 + 1]$. For the rest of the proof to be correct, one *needs* $\varepsilon(r) \geq \varepsilon$, which is certainly true, provided $\varepsilon(r) \downarrow 0$ sufficiently slowly.

As defining such a sequence $\varepsilon(r) \downarrow 0$ may be difficult, a weaker result can be obtained for an arbitrary (but fixed) precision $\varepsilon_0$, by fixing *a priori* $\varepsilon(r) = \varepsilon$ for all $r$. Indeed, provided $\varepsilon$ is sufficiently small, one then obtains $|\inf \mathbb{Q}_r - \rho^*| < \varepsilon_0$ for all $r$ sufficiently large.

Finally, on page 1131, line 15 from bottom, "an admissible" should be "a strictly admissible;" on page 1135, lines 8–9, interchange "$\mathbb{D}$" with "$\mathbb{P}$."

## REFERENCE

Lasserre, J. B. (2002). Bounds on measures satisfying moment conditions. *Ann. Appl. Probab.* **12** 1114–1137. MR1925454

LAAS-CNRS
7 Avenue du Colonel Roche
31077 Toulouse Cédex 4
France
e-mail: lasserre@laas.fr